\documentclass[reqno,10pt]{amsproc}

\usepackage{amssymb}
\usepackage{graphicx}
\usepackage{amscd}
\usepackage{amsmath}

\theoremstyle{plain}

\newtheorem{convention}{Convention}

\newtheorem{lemma}{Lemma}

\newtheorem{theorem}{Theorem}
\numberwithin{equation}{section}

\begin{document}

\title{Unified approach to classical equations of inverse problem theory}
\author{ M. I. Belishev}
\address{St.Petersburg   Department   of   V.A.Steklov    Institute   of   Mathematics
of   the   Russian   Academy   of   Sciences, 7, Fontanka, 191023
St.Petersburg, Russia} \email{belishev@pdmi.ras.ru, }
\author{ V. S. Mikhaylov}
\address{St.Petersburg   Department   of   V.A.Steklov    Institute   of   Mathematics
of   the   Russian   Academy   of   Sciences, 7, Fontanka, 191023
St.Petersburg, Russia and Chebyshev Laboratory, St. Petersburg
State University, 14th Line, 29b, Saint Petersburg, 199178 Russia}
\email{vsmikhaylov@pdmi.ras.ru, }
\thanks{the work of the authors was supported by RFFI 11-01-00407A, NSh-4210.2010.1,
and the Chebyshev Laboratory (Department of
Mathematics and Mechanics, St. Petersburg State University) under
RF Government grant 11.G34.31.0026}
\date{May 17, 2012}

\maketitle

\noindent {\bf Abstract.} The boundary control (BC-) method is an
approach to inverse problems based upon their deep relations to
control and system theory. We show that the classical integral
equations of inverse problem theory (Gelfand-Levitan, Krein and
Marchenko equations) can be derived in the framework of the
BC-method in a unified way. Namely, to solve each of these
equations is in fact to solve a relevant boundary control problem,
whereas its solution is determined by the inverse data.

\section{Introduction.}

Our paper is mainly of methodical character. Its aim is to
demonstrate that the well-known classical integral equations of
inverse problem theory can be derived in the framework of a
unified approach. The approach is the boundary control (BC-)
method that is a method for solving inverse problems based on
their relations to control and system theory \cite{B97, B07, B08}.
As we show, to solve the classical equations is in fact to solve
boundary control problems (BCP) of a special kind.

The latter problems are set up for dynamical systems governed by
the wave equation $u_{tt}+Lu=0$ on the semi-axis $x>0$ with a
boundary control $f$ at $x=0$ of Dirichlet or Neumann type. Here
$L$ is a Sturm-Liouville type operator. The solution $u=u^f(x,t)$
describes a wave produced by the control. A BCP is to find a
control, which transfers the system from the zero initial state to
a given final state $y$, i.e., provides $u^f(\,\cdot\,,T)=y$. A
{\it special BCP} corresponds to the case of $y$ provided $Ly=0$.
Namely,

\noindent $\bullet$\,\, 
Let $L=-\frac{d^2}{dx^2}+q(x)$, and $y$ satisfy $Ly=0$,
$y(0)=0,\,\,y^\prime(0)=1$. We show that the Dirichlet boundary
control $u|_{x=0}=f(t)$, which solves the BCP
$u^f(\,\cdot\,,T)=y$, coincides (up to a simple explicit
transformation) with the solution of the Gelfand-Levitan equation
derived in \cite{GL} for solving the inverse spectral
Sturm-Liouville problem.

\noindent $\bullet$\,\, 
Let $L=-\rho^{-1}(x)\frac{d^2}{dx^2}$, and $y$ satisfy $Ly=0,\,
y(0)=1,\,\,y^\prime(0)=0$. We show that the Neumann boundary
control $u_x|_{x=0}=f(t)$, which solves the BCP
$u^f_t(\,\cdot\,,T)=y$, coincides (up to a trivial change of
variables) with the solution of the M.Krein equation
\cite{Kr1,Kr2}, which is used for solving the inverse spectral
problem for inhomogeneous string. Note that at this point, in
fact, we reproduce the results of Sondhi and Gopinath
\cite{GoSh1}, \cite{GoSh2}.

Also, taking $Ly=0,\,\,\,y(0)=0,\,\,y^\prime(0)=0$, we derive the
so-called Pariiskii equation \cite{Bl71}.

\noindent $\bullet$\,\, At last, we consider the dynamical
scattering problem for the equation $u_{tt}-u_{xx}+qu=0$ on the
semi-axis $x>0$. It is shown that the solution of the Marchenko
equation also solves a relevant BCP, in which the role of controls
is played by the incoming waves.

Such a unified look at the classical equations is quite
productive: it provides their natural analogs in multidimensional
problems \cite{B87} and relates the BC-method with another
approaches \cite{AM}.

\section{Forward problems}

\subsection{Problem $1_{\rm Dir}$}
We deal with an initial boundary value problem $1_{\rm Dir}$ of
the form
\begin{align}
& \rho u_{tt}-u_{xx}+qu=0, && x>0,\,\, t>0\label{eq_1} \\
& u|_{t=0}=u_t|_{t=0}=0, && x\geqslant 0 \label{ic_1}\\
& u|_{x=0}=f, && t\geqslant 0, \label{bc_1}
\end{align}
where $\rho\in C^\infty[0,+\infty)$, $\rho(x)>0$  is a {\it
density}; $q\in C^\infty[0,+\infty)$ is a {\it potential}
\footnote{In the paper, all functions real. Everywhere {\it
smooth} means $C^\infty-${\it smooth}.} depending on $x$; $f=f(t)$
is a Dirichlet boundary control; $u=u^f(x,t)$ is a solution ({\it
wave}). Also, we assume that
\begin{align}\label{Eikonal infty}
\int_0^\infty\rho^{\frac{1}{2}}(x)\,dx=\infty.
\end{align}
Under this assumption, problem (\ref{eq_1})--(\ref{bc_1}) can be
reduced to the relevant Volterra-type integral equations and turns
out to be well posed (see, e.g., \cite{B08, Blag2}). Then, one can
show that for a smooth $f$ vanishing near $t=0$ the problem has a
unique classical smooth solution $u=u^f(x,t)$.

The same integral equations technique allows to introduce a
fundamental solution $u^\delta(x,t)$ for $f=\delta(t)$ (the Dirac
delta-function) and derive the well-known representation
\begin{equation*}
u^\delta(x,t)=\left[\frac{\rho(x)}{\rho(0)}\right]^{-\frac{1}{4}}\delta(t-\tau(x))+w(x,t),
\end{equation*}
where $\tau(x):=\int_0^x\rho^{\frac{1}{2}}(s)\,ds$ is an
\emph{eikonal} satisfying $\tau(\infty)=\infty,$ the function $w$
is smooth in $\{(x,t)\,|\, 0\leqslant \tau(x)\leqslant t\}$ and
vanishes as $t>\tau(x).$

For the classical solution, the Duhamel representation holds:
\begin{align}
& u^f(x,t)=u^\delta(x,t)*f(t)=\notag\\
\label{u_Duhamel}&
\left[\frac{\rho(x)}{\rho(0)}\right]^{-\frac{1}{4}}f(t-\tau(x))+\int_0^{t-\tau(x)}w(x,t-s)f(s)\,ds,\quad
x\geqslant 0,\,\,t\geqslant 0.
\end{align}
\begin{convention}
In $(\ref{u_Duhamel})$ and everywhere in the sections 1--4, we
assume that all functions depending on time are extended to $t<0$
by zero.
\end{convention}
\noindent For the case $f\in L_2^{\rm loc}(0,\infty)$, which we
mainly deal with, the (generalized) solution $u=u^f(x,t)$ of the
problem $1_{\rm Dir}$ {\it is defined} as the right hand side of
(\ref{u_Duhamel}). As is easy to see, such a solution belongs to
the class $C_{\rm loc}\left([0,\infty),L_2(0,\infty)\right)$ (as
an $L_2(0,\infty)-$valued function of $t$) and satisfies
$u^f|_{t<\tau(x)}=0$.

By the latter, we have
\begin{equation}
\label{u_supp} \operatorname{supp}u^f(\,\cdot\,,t)\subset
[0,x(t)],\qquad t\geqslant 0,
\end{equation}
where $x(\tau)$ is the inverse function of the eikonal $\tau(x)$.
This means that the waves $u^f$ propagate along the semi-axis
$x\geqslant 0$ with a finite (variable) velocity
$c=\frac{dx}{dt}|_{t=\tau(x)}=\frac{d\tau}{dx}=\rho^{-\frac{1}{2}}(x).$
Owing to assumption (\ref{Eikonal infty}), the waves
$u^f(\,\cdot\,,t)$ are compactly supported for any $t\geqslant 0$.

For an $s>0$, define a {\it delay operator} ${\mathcal T}_s$,
which acts on the time variable by $({\mathcal T}_s h)(t)=h(t-s)$
(recall Convention 1!). Let $J:= \int^t_0$ be integration in time.
An independence of the operator $L=\rho^{-1}(-\frac{d^2}{dx^2}+q)$
from time leads to the well-known properties of solutions:
\begin{equation}
\label{u diff} u^{{\mathcal T}_s f}={\mathcal T}_s u^f, \quad
u^{Jf}=Ju^f, \quad u^{\frac{df}{dt}}=u^f_{t}, \quad
u^{\frac{d^2f}{dt^2}}=u^f_{tt}=Lu^f\,.
\end{equation}
Note that the second, third, and forth relations follow from the
first one.

\subsection{Problems $1^T_{\rm Dir}$ and $\widetilde 1^{2T}_{\rm Dir}$}
The problem
\begin{align}
& \rho u_{tt}-u_{xx}+qu=0, &&  0<x<x(T),\,\, 0<t<T\notag\\ 
& u|_{t<\tau(x)}=0,\notag\\ 
& u|_{x=0}=f, && 0\leqslant t \leqslant T \notag 
\end{align}
with a final moment $t=T<\infty$ is referred to as a problem
$1^T_{\rm Dir}$. Under the condition (\ref{Eikonal infty}), it is
well posed for any $T>0$ and meaningful for $f\in L_2(0,T).$

By the same reason, the problem $\widetilde 1^{2T}_{\rm Dir}$ of
the form
\begin{align}
& \rho u_{tt}-u_{xx}+qu=0, && 0<x<x(T),\,\, 0<t<2T-\tau(x)\label{eq_1_2T} \\
& u|_{t<\tau(x)}=0 \label{ic_1_2T}\\
& u|_{x=0}=f, && 0\leqslant t \leqslant 2T, \label{bc_1_2T}
\end{align}
is also well-posed. It can be regarded as a natural extension of
the problem $1^T_{\rm Dir}$, which exists owing to the finiteness
of the wave propagation speed \footnote{such an extension is
standard in hyperbolic problems \cite{Blag2}}. Its solution $u^f$
is also represented by (\ref{u_Duhamel}) (for the relevant values
of $x$ and $t$).

An evident but important fact is that the solutions of the
problems $1^T_{\rm Dir}$ and $\widetilde 1^{2T}_{\rm Dir}$ are
determined by the values of the coefficients $\rho,$ $q$ on the
segment $[0,x(T)]$ only (do not depend on $\rho,q|_{x>x(T)}$).

\subsection{Problems with Neumann control}
Under the same assumptions on $\rho$ and $q$, we consider a
problem $1_{\rm Neum}$ of the form
\begin{align}
& \rho u_{tt}-u_{xx}+qu=0, &&  x>0,\,\, t>0 \notag\\ 
& u|_{t=0}=u_t|_{t=0}=0, && x\geqslant 0\notag \\ 
& u_{x}|_{x=0}=f, && t\geqslant 0 \notag 
\end{align}
with a Neumann boundary control $f$. This problem can also be
reduced to the integral Volterra-type equations. For
$f=\delta(t)$, its fundamental solution is of the well-known form
\begin{equation}
\label{u_2_delta}
u^\delta(x,t)=-\left[\rho(0)\rho(x)\right]^{-\frac{1}{4}}\theta(t-\tau(x))+w(x,t),
\end{equation}
where $\theta$ is the Heavyside function: $\theta|_{t<0}=0,$
$\theta|_{t\geqslant 0}=1,$ $w$ is smooth in $\{(x,t)\,|\,
0\leqslant \tau(x)\leqslant t\}$, continuous in $\{(x,t)\,|\,
x\geqslant 0,\, t\geqslant 0\}$ and vanishes as $t \leqslant
\tau(x).$

For $f\in L_2^{\rm loc}(0,\infty)$, we {\it define} a solution by
Duhamel
\begin{equation}
\label{u_2_Duhamel} u^f(x,t):=u^\delta(x,t)*f(t)=\langle {\rm
see}\,
 (\ref{u_2_delta})\rangle=\int_0^{t-\tau(x)}r(x,t-s)f(s)\,ds
\end{equation}
for all $x\geqslant 0,\,\,t\geqslant 0$, where $r:=u^\delta$ is
smooth in $\{(x,t)\,|\, 0\leqslant \tau(x)\leqslant t\}$, vanishes
as $t< \tau(x)$, and satisfies
\begin{equation}\label{jump of r}
r(x,\tau(x)+0)= \langle{\rm see\,}(\ref{u_2_delta})\rangle =
-\left[\rho(0)\rho(x)\right]^{-\frac{1}{4}}\not=0\,.
\end{equation}
For such a solution, we have
\begin{equation*}
\operatorname{supp}u^f(\,\cdot\,,t)\subset
[0,x(t)],\qquad t\geqslant 0.
\end{equation*}
Also, $u^f(\,\cdot\,,t)$ belongs to
$H^1[0,\infty)$\footnote{$H^1[0,\infty):=\{y\,|\,y,\,y^\prime \in
L_2(0,\infty)\}$ is the Sobolev class.} for any $t \geqslant 0$
and satisfies $u^f|_{t\leqslant \tau(x)}=0$.
\smallskip

Owing to the finiteness of the wave propagation speed, the problem
$1^T_{\rm Neum}$
\begin{align}
& \rho u_{tt}-u_{xx}+qu=0, &&  0<x<x(T),\,\, 0<t<T\notag \\ 
& u|_{t<\tau(x)}=0 \notag \\ 
& u_x|_{x=0}=f, && 0\leqslant t \leqslant T \notag 
\end{align}
and its extension by hyperbolicity $\widetilde 1^{2T}_{\rm Neum}$
\begin{align}
& \rho u_{tt}-u_{xx}+qu=0, && 0<x<x(T),\,\, 0<t<2T-\tau(x)\notag \\ 
& u|_{t<\tau(x)}=0 \notag\\ 
& u_x|_{x=0}=f, && 0\leqslant t \leqslant 2T \notag 
\end{align}
turn out to be well-posed for the above introduced class of
controls. Their solutions are determined by $\rho,q|_{[0,x(T)]}$
(do not depend on $\rho,q|_{x>x(t)}$).

\section{Boundary control problems}

\subsection{System $\alpha_{\rm Dir}^T$}
Here we consider the problem $1^T_{\rm Dir}$ as a dynamical system
(that we call $\alpha_{\rm Dir}^T$) and endow it with the standard
attributes of control theory: spaces and operators.

\noindent {\bf Outer space}\,\,\,A Hilbert space of controls
$\mathcal{F}^T:=L_2(0,T)$ with the inner product
\begin{equation*}
(f,g)_{\mathcal{F}^T}=\int_0^T f(t)g(t)\,dt
\end{equation*}
is said to be an \emph{outer space} of the system $\alpha_{\rm
Dir}^T.$ It contains a family of subspaces
\begin{equation*}
{\mathcal F}^{T,\,\xi}:=\left\{f\in {\mathcal F}^T\,|\, {\rm supp}
f\subset [T-\xi,T]\right\},\quad 0\leqslant \xi\leqslant T,
\end{equation*}
which consists of the delayed controls. Here $T-\xi$ is the delay,
$\xi$ is the action time. Note that ${\mathcal F}^{T,\,0}=\{0\}$
and ${\mathcal F}^{T,\,T}={\mathcal F}^T$.
\smallskip

\noindent {\bf Inner space}\,\,\, A wave $u^f(\,\cdot\,,t)$ is
interpreted as a {\it state} of the system at the moment $t$. A
space $\mathcal{H}^T:=L_{2,\,\rho}(0,x(T))$ with the product
\begin{equation*}
(u,v)_{\mathcal{H}^T}=\int_0^{x(T)} u(x)v(x)\,\rho(x)\,dx
\end{equation*}
is an \emph{inner space}. By $L_2-$regularity of solutions to
$1^T_{\rm Dir}$ and (\ref{u_supp}), one has $u^f(\cdot,t)\in
\mathcal{H}^T$ for all $t\in [0,T]$. The inner space contains a
family of subspaces
\begin{equation*}
\mathcal{H}^{\xi}:=\left\{y\in \mathcal{H}^T\,|\, {\rm supp\,}
y\subset [0,x(\xi)]\right\},\quad 0\leqslant \xi\leqslant T.
\end{equation*}
By (\ref{u_supp}) and the first relation in (\ref{u diff}), $f\in
\mathcal{F}^{T,\,\xi}$ implies $u^f(\,\cdot\,,T)\in
\mathcal{H}^{\xi}$.
\smallskip

\noindent {\bf Control operator}\,\,\,In the system $\alpha_{\rm
Dir}^T$, the `input $\mapsto$ state' correspondence is realized by
a \emph{control operator} $W^T:\mathcal{F}^T\to \mathcal{H}^T,\,\,
W^Tf:=u^f(\,\cdot\,,T)$. With regard to (\ref{u_Duhamel}), we have
the representation
\begin{equation}
\label{W_T_repr}\left(W^T
f\right)(x)=\left[\frac{\rho(x)}{\rho(0)}\right]^{-\frac{1}{4}}f(T-\tau(x))+
\int_0^{T-\tau(x)}w(x,T-s)f(s)\,ds, \quad 0\leqslant x \leqslant
x(T).
\end{equation}
\begin{lemma}
\label{W_T_lemma} The operator $W^T$ is an isomorphism from
$\mathcal{F}^T$ onto $\mathcal{H}^T$. The relation
$W^T\mathcal{F}^{T,\,\xi}=\mathcal{H}^\xi$ is valid as
$0\leqslant\xi\leqslant T$.
\end{lemma}
\begin{proof}\,\,\,The continuity of $W^T$ is evident from
(\ref{W_T_repr}). Choose a $y\in \mathcal{H}^T$ and consider the
equation $W^Tf=y$. Substituting $t=T-\tau(x)$ in (\ref{W_T_repr}),
we have
\begin{equation*}
\alpha(t)f(t)+\int_0^t \widetilde w(x,s)f(s)\,ds=y(x(T-t)), \quad
0\leqslant t\leqslant T
\end{equation*}
with a strictly positive $\alpha$ and smooth $\widetilde w$, i.e.,
get a second order integral Volterra equation. The latter is
uniquely solvable and, hence, $(W^T)^{-1}$ does exist and is
continuous.

Since $T>0$ is arbitrary, for each system $\alpha_{\rm Dir}^\xi$
one has $W^\xi\mathcal{F}^\xi=\mathcal{H}^\xi$ as $\xi \leqslant
T$. By the delay relation in (\ref{u diff}), this is equivalent to
$W^T\mathcal{F}^{T,\,\xi}=\mathcal{H}^\xi$.
\end{proof}
\smallskip

\noindent {\bf Response operators}\,\,\, The `input $\mapsto$
output' correspondence is realized by a \emph{response operator}
$R^T: \mathcal{F}^T\to \mathcal{F}^T,\,\,{\rm
Dom\,}R^T=\left\{f\in H^1[0,T]\,|\, f(0)=0\right\}$, \,\,$R^Tf
:=u^f_x(0,\,\cdot\,)$. By the use of (\ref{u_Duhamel}), one can
represent
\begin{equation}
\label{R_T_repr} \left(R^T
f\right)(t)=-\alpha\frac{df}{dt}(t)+\beta f(t)+\int_0^t
r(t-s)f(s)\,ds,\qquad 0\leqslant t\leqslant T,
\end{equation}
where $\alpha=\rho^{\frac{1}{2}}(0),$
$\beta=-[4\rho(0)]^{-1}\rho^\prime(0),$ $r=w_x(0, \,\cdot\,)\in
C^\infty[0,T].$

Also, the system $\alpha_{\rm Dir}^T$ is endowed with an
\emph{extended response operator}. It is associated with the
problem $\widetilde 1^{2T}_{\rm Dir}$ and introduced by
$\widetilde R^{2T}: \mathcal{F}^{2T}\to \mathcal{F}^{2T}$, ${\rm
Dom\,}\widetilde R^{2T}=\left\{f\in H^1[0,2T]\,|\,
f(0)=0\right\}$, $\widetilde R^{2T}f:=u^f_x(0,\,\cdot\,)$, where
$u^f$ is the solution to (\ref{eq_1_2T})--(\ref{bc_1_2T}). Since
(\ref{u_Duhamel}) holds for the extended problem, we have
\begin{equation}
\label{R_2T_repr} (\widetilde
R^{2T}f)(t)=-\alpha\frac{df}{dt}(t)+\beta f(t)+\int_0^t
r(t-s)f(s)\,ds,\qquad 0\leqslant t\leqslant 2T
\end{equation}
with the same $\alpha,$ $\beta$ and $r\in
C^\infty[0,2T]$\footnote{as is easy to recognize, reducing this
$r$ on the subinterval $[0,T]$ one gets the kernel $r$ from
(\ref{R_T_repr})}. Along with the solution $u^f$ of the problem
$\widetilde 1^{2T}_{\rm Dir}$, the operator $\widetilde R^{2T}$ is
determined by the values of the coefficients $\rho, r|_{[0,x(T)]}$
and does not depend on $\rho, r|_{x>x(T)}$. That is why it is
natural to regard $\widetilde R^{2T}$ as an intrinsic object of
the system $\alpha^T_{\rm Dir}$ (but not $\alpha^{2T}_{\rm Dir}$).

The equalities (\ref{u diff}) easily lead to the commutation
relations \begin{equation}\label{commut} R {\mathcal T}_s =
{\mathcal T}_s R, \qquad R J = J R\,,
\end{equation}
which are valid for $R=R^T$ and $R=\widetilde R^{2T}$ on the
relevant time intervals. Also, it is clear that the knowledge of
$\widetilde R^{2T}$ suffices to know $R^\xi$ and $\widetilde
R^{2\xi}$ for every $\xi \leqslant T$.
\smallskip

\noindent {\bf Connecting operator}\,\,\, A map $C^T:
\mathcal{F}^T\to \mathcal{F}^T,\,\,  C^T:=(W^T)^*W^T$ is said to
be a \emph{connecting operator} of the system $\alpha^T_{\rm
Dir}.$ By this definition, for $f,g\in \mathcal{F}^T$ one has
\begin{equation}
\label{C_T}
\left(C^Tf,g\right)_{\mathcal{F}^T}=\left(W^Tf,W^Tg\right)_{\mathcal{H}^T}=
\left(u^f(\,\cdot\,,T),u^g(\,\cdot\,,T)\right)_{\mathcal{H}^T},
\end{equation}
i.e., this operator connects the metrics of the outer and inner
spaces. Since $W^T$ is an isomorphism from $\mathcal{F}^T$ onto
$\mathcal{H}^T$ (Lemma \ref{W_T_lemma}), the connecting operator
$C^T=(W^T)^*W^T$ turns out to be a positive definite isomorphism
of $\mathcal{F}^T$, i.e., satisfies
$(C^Tf,f)_{\mathcal{F}^T}\geqslant c\|f\|^2$ with a constant
$c>0$.

A basic fact of the BC-method is a simple and explicit relation
between the operators $C^T$ and $\widetilde R^{2T}$. Introduce the
auxiliary operators $S^T: \mathcal{F}^T\to \mathcal{F}^{2T}$,
\begin{equation*}
\left(S^Tf\right)(t):=\left\{\begin{array}l \quad f(t),\quad
\qquad 0\leqslant
t\leqslant T\\
-f(2T-t),\quad T < t\leqslant 2T
\end{array}
\right.
\end{equation*}
and $J^{2T}:\mathcal{F}^{2T}\to
\mathcal{F}^{2T},\,\,(J^{2T}f)(t):=\int_0^tf(s)\,ds,\quad
0\leqslant t\leqslant 2T$. For the operator
$(S^T)^*:\mathcal{F}^{2T}\to \mathcal{F}^{T}$, the equality
$(S^T)^*=2N^{2T} P_-^{2T}$ is valid, where
$P_-^{2T}:\mathcal{F}^{2T}\to \mathcal{F}^{2T}$ projects on the
subspace of odd (w.r.t. $t=T$) functions, and
$N^{2T}:\mathcal{F}^{2T}\to \mathcal{F}^{T}$ reduces functions
from $[0,2T]$ to $[0,T]$.
\begin{theorem}
\label{C_T_theor} The relation
\begin{equation}
\label{C_T_repr} C^T=\frac{1}{2}(S^T)^*\widetilde R^{2T}J^{2T}S^T
\end{equation}
and representation
\begin{equation}
\label{C_T_repr1} \left(C^Tf\right)(t)=\alpha
f(t)+\int_0^T\left[\frac{1}{2}\int_{|t-s|}^{2T-t-s}r(\eta)\,d\eta\right]f(s)\,ds
\end{equation}
are valid, where $r$ is taken from $(\ref{R_2T_repr})$.
\end{theorem}
\begin{proof}\,\, Choose arbitrary controls $f,g\in
C_0^\infty(0,T)$ and denote $f_-:=S^Tf\in C_0^\infty(0,2T)$. Let
$u^{f_-}$ be the solution to the problem $\widetilde 1^{2T}_{\rm
Dir}$ and $u^g$ the solution to $1^T_{\rm Dir}$. Note that by the
choice of $g$ and (\ref{u_supp}) the latter solution satisfies
$\operatorname{supp} u^g(\cdot,t)\subset[0,x(T))$ (i.e. vanishes
near $x=x(T)$) for all $t\in [0,T]$.

It is easy to see that the \emph{Blagoveschenskii function}
\begin{equation*}
b(s,t)=\int_0^{x(T)}u^{f_-}(x,s)u^g(x,t)\,\rho(x)\,dx
\end{equation*}
is well defined in the triangular $\left\{(s,t)\,|\, 0\leqslant
t\leqslant T,\,\, t\leqslant s\leqslant 2T-t\right\}$ and
satisfies
\begin{align*}
& b_{tt}(s,t)-b_{ss}(s,t)=\int_0^{x(T)}\left[u^{f_-}(x,s)u_{tt}^g(x,t)-
u_{ss}^{f_-}(x,s)u^g(x,t)\right]\rho(x)\,dx=\notag\\
& \int_0^{x(T)}\left\{u^{f_-}(x,s)\left[u_{xx}^g(x,t)-q(x)u^g(x,t)\right]-
\left[u_{xx}^{f_-}(x,s)-q(x)u^{f_-}(x,s)\right]u^g(x,t)\right\}\,dx=\notag\\
&-u^{f_-}(0,s)u_x^g(0,t)+u_x^{f_-}(0,s)u^g(0,t)=-f_-(s)(R^Tg)(t)+(\widetilde
R^{2T}f_-)(s)g(t)=:h(s,t)\,.\label{b_1}
\end{align*}
Taking into account $b(s,0)=b_t(s,0)=0$ for $0\leqslant s\leqslant
2T,$ and integrating the equation $b_{tt}-b_{ss}=h$ by D'Alembert,
with regard to the oddness of $f_-,$ we get
\begin{align*}
&
b(T,T)=\frac{1}{2}\int_0^T\,dt\,\int_t^{2T-t}\left[-f_-(s)(R^Tg)(t)+(\widetilde
R^{2T}f_-)(s)g(t)\right]\,ds\\
&
=\int_0^T\,dt\,g(t)\left[-\frac{1}{2}\left(\int_0^t-\int_0^{2T-t}\right)(\widetilde
R^{2T}f_-)(s)\,ds\right]=\left(g,
\frac{1}{2}\left(S^T\right)^*J^{2T}\widetilde
R^{2T}S^Tf\right)_{\mathcal{F}^T}.
\end{align*}

The solutions $u^f$ and $u^{f_-}$ of the problems $1^T_{\rm Dir}$
and $\widetilde 1^{2T}_{\rm Dir}$ respectively coincide as
$0\leqslant t\leqslant T.$ Hence, we have
\begin{align*}
& b(T,T)=\int_0^{x(T)}u^{f_-}(x,T)u^g(x,T)\rho(x)\,dx=\int_0^{x(T)}u^{f}(x,T)u^g(x,T)\rho(x)\,dx\\
&=\left(u^f(\,\cdot\,,T),u^g(\,\cdot\,,T)\right)_{\mathcal{H}^T}=\langle
{\rm see\,}
(\ref{C_T})\rangle=\left(C^Tf,g\right)_{\mathcal{F}^T}.
\end{align*}
Comparing two expressions obtained for $b(T,T)$, with regard to a
density of $C_0^\infty(0,T)$ in $\mathcal{F}^T$, we get
\begin{equation*}
C^Tf=\frac{1}{2}\,(S^T)^*J^{2T}\widetilde R^{2T}S^T, \quad f\in
C_0^\infty(0,T).
\end{equation*}
By (\ref{commut}), we have $J^{2T}\widetilde R^{2T}=\widetilde
R^{2T}J^{2T}$. As is evident from (\ref{R_2T_repr}), $\widetilde
R^{2T}J^{2T}$ is a continuous operator defined on the whole
$\mathcal{F}^T$ (not only on $\operatorname{Dom}\widetilde
R^{2T}$). Hence, we arrive at (\ref{C_T_repr}). Substituting
(\ref{R_2T_repr}) into (\ref{C_T_repr}), after simple calculations
we get (\ref{C_T_repr1}).
\end{proof}

\subsection{BCP for $\alpha^T_{\rm Dir}.$}
A boundary control problem (BCP) for the system $\alpha^T_{\rm
Dir}$ is:

\noindent $\bullet$\,\,\,{\it given a target state $y\in
\mathcal{H}^T$, find a control $f\in \mathcal{F}^T$ such that
$u^f(\,\cdot\,,T)=y$}.

\noindent In other words, the question is whether one can manage a
shape of the wave by managing the boundary control.

Such a problem is equivalent to the equation $W^Tf=y.$ As a
straightforward consequence of Lemma \ref{W_T_lemma}, we conclude
that the BCP is well posed: it has a unique solution
$f=(W^T)^{-1}y$ for any $y\in \mathcal{H}^T$, the solution
depending on $y$ continuously (in the relevant $L_2-$norms).
Moreover, if $y\in {\mathcal H}^\xi$ then one has
$f=(W^T)^{-1}y\in \mathcal{F}^{T,\,\xi}$ ($0\leqslant \xi\leqslant
T$).

Now, consider a \emph{special} BCP. Assume that the target
function $y$ is a solution of the Sturm-Liouville equation
$-y^{\prime \prime}+qy=0$ as $x\in(0,x(T))$.
\begin{theorem}
\label{Theor_DIR_BCP} The control $f=\left(W^T\right)^{-1}y$ is a
unique solution to the equation
\begin{equation}
\label{C_T_eqn} C^Tf=y^\prime(0)\varkappa^T-y(0)(R^T)^*\varkappa^T
\end{equation}
in $\mathcal{F}^T$, where $\varkappa^T(t):=T-t.$
\end{theorem}
\begin{proof} Choose a $g\in C_0^\infty(0,T)$; let $u^g$
be the solution to the problem $1^T_{\rm Dir}$. Since $g$ vanishes
near $t=0$ one has $\operatorname{supp}u^g(\,\cdot\,,t)\subset
[0,x(T))$, i.e., the wave $u^g(\,\cdot\,,t)$ vanishes near
$x=x(T)$ for all $t\in[0,T]$. For $f=(W^T)^{-1}y$, one has
\begin{align*}
&(C^Tf,g)_{\mathcal{F}^T}=\langle {\rm see\,}
(\ref{C_T})\rangle=(W^Tf,W^Tg)_{\mathcal{H}^T}=(y,u^g(\,\cdot\,,T))_{\mathcal{H}^T}=\\
&\int_0^{x(T)} y(x)u^g(x,T)\rho(x)\,dx=\int_0^{x(T)}
y(x)\left[\int_0^T(T-t)u^g_{tt}(x,t)\,dt\right]\rho(x)\,dx=\\
&\int_0^T\,dt\,
\varkappa^T(t)\int_0^{x(T)}\rho(x)u_{tt}^g(x,t)y(x)\,dx=\\
&\int_0^T\,dt \,
\varkappa^T(t)\int_0^{x(T)}\left[u_{xx}^g(x,t)-q(x)u^g(x,t)\right]y(x)\,dx=\\
&\int_0^T\,dt \,
\varkappa^T(t)\left\{-u^g_x(0,t)y(0)+u^g(0,t)y'(0)-\int_0^{x(T)}u^g(x,t)\left[y''(x)-q(x)y(x)\right]\,dx\right\}\\
&= \left(\varkappa^T, -y(0)R^T g+ y'(0)g\right)_{{\mathcal
F}^T}=\left(y'(0)\varkappa^T-y(0)(R^T)^*\varkappa^T,g\right)_{\mathcal{F}^T}\,.
\end{align*}
By arbitrariness of $g$ we get (\ref{C_T_eqn}).
\end{proof}

Taking into account (\ref{C_T_repr1}), we see that (\ref{C_T_eqn})
is a second order Fredholm integral equation, which is uniquely
solvable by the properties of  $C^T$. Also, a simple analysis
shows that the solution $f$ is of the class $C^\infty[0,T].$

For the solution of the special BCP, the relation
(\ref{u_Duhamel}) takes the form
\begin{equation*}
y(x)=u^f(x,T)=\left[\frac{\rho(x)}{\rho(0)}\right]^{-\frac{1}{4}}f(T-\tau(x))+\int_0^{T-\tau(x)}w(x,T-s)f(s)\,ds.
\end{equation*}
Putting here $x=x(T)$ (so that $\tau(x)=T$), we obtain the
relation
\begin{equation}
\label{Jump_cnd}
y(x(T))=\left[\frac{\rho(x(T))}{\rho(0)}\right]^{-\frac{1}{4}}f(0),
\end{equation}
which will be used for solving the inverse problems.

In fact, in inverse problems, not a single equation but a family
of equations of the form (\ref{C_T_eqn}) is in use. Namely,
substituting $T$ for an intermediate $\xi\in [0,T]$, we have the
equation
\begin{equation}
\label{C_T_del_eqn} C^\xi
f=y^\prime(0)\varkappa^\xi-y(0)(R^\xi)^*\varkappa^\xi
\end{equation}
in $\mathcal{F}^\xi$, where $\varkappa^\xi(t)=\xi-t.$ For their
solutions $f=f^\xi(t)$, relation (\ref{Jump_cnd}) takes the form
\begin{equation}
\label{Jump_cnd_xi}
y(x(\xi))=\left[\frac{\rho(x(\xi))}{\rho(0)}\right]^{-\frac{1}{4}}f^\xi(0),
\qquad 0\leqslant \xi \leqslant T.
\end{equation}
As we mentioned in sec 3.1, the operator $\widetilde R^{2T}$
determines the operators $\widetilde R^{2\xi}$ and $R^\xi$ (and,
hence, $C^\xi$) for all $\xi<T$. Therefore, the knowledge of
$\widetilde R^{2T}$ suffices to use equations (\ref{C_T_del_eqn})
and their solutions $f^\xi.$

The solutions $f^\xi$ solve the BCPs
\begin{equation*}
u^f(x,\xi)=(W^\xi f)(x)=y(x),\quad 0\leqslant x\leqslant x(\xi).
\end{equation*}
For what follows, it it convenient to refer them to the mutual
final moment $t=T$: by the first relation in (\ref{u diff}), the
controls $f^{T,\,\xi}\in \mathcal{F}^{T,\,\xi}\subset
\mathcal{F}^T$ of the form \footnote{recall Convention 1!}
\begin{equation*}
f^{T,\,\xi}(t):=f^\xi(t-(T-\xi)),\quad 0 \leqslant
t \leqslant T
\end{equation*}
solve the BCPs
\begin{equation}
\label{BCP_del} u^f(x,T)=(W^Tf)(x)=y^\xi(x):=\left\{\begin{array}l
y(x),\qquad 0\leqslant
x\leqslant x(\xi)\\
0,\qquad x(\xi)< x\leqslant x(T).
\end{array}
\right.
\end{equation}

\subsection{System $\alpha^T_{\rm Neum}.$}
A dynamical system with Neumann boundary control is associated
with the problem $1^T_{\rm Neum}$. Its {\bf outer} and {\bf inner}
spaces are introduced in the same way: $\mathcal{F}^T=L_2(0,T)$,
${\mathcal H}^T=L_{2,\,\rho}(0,x(T)).$

By (\ref{u_2_Duhamel}), a {\bf control operator} $W^T:
\mathcal{F}^T\to \mathcal{H}^T$, $W^Tf:=u^f(\,\cdot\,,T)$ is of
the form
\begin{equation}
\label{W_T_repr_NM}
\left(W^Tf\right)(x)=\int_0^{T-\tau(x)}r(x,T-s)f(s)\,ds,\quad
0\leqslant x\leqslant x(T).
\end{equation}
Thus, $W^T$ is a compact operator. In the mean time, for an
$f\in{\mathcal F}^T_+:=\{f \in H^1[0,T]\,|\,\,f(0)=0\}$,
integration by parts leads to
\begin{equation}\label{repr Wdf/dt}
\left(W^T\frac{df}{dt}\right)(x)=
\alpha(x)f(T-\tau(x))+\int_0^{T-\tau(x)}r_t(x,T-s)f(s)\,ds
\end{equation}
with $\alpha(x):=-r(x,\tau(x)+0)\not=0$ (see (\ref{jump of r})).
We omit the proof of the following facts, which can be easily
derived from representations (\ref{W_T_repr_NM}) and (\ref{repr
Wdf/dt}).
\begin{lemma}
\label{Lemma_NM_contr} The operator $W^T$ maps ${\mathcal F}^T$
onto $\left\{y\in H^1[0,x(T)]\,|\,\,y(x(T))=0\right\}$
injectively. The operator $W^T\!\frac{d}{dt}$ defined on
${\mathcal F}^T_+$ is continuous. Its extension by continuity maps
$\mathcal{F}^T$ onto $\mathcal{H}^T$ isomorphically, the relations
$W^T\!\frac{d}{dt} {\mathcal F}^{T,\,\xi}={\mathcal H}^\xi$ being
valid for all $0\leqslant \xi \leqslant T$.
\end{lemma}
\smallskip

A {\bf response operator} of the system $\alpha^T_{\rm Neum}$ is
$R^T: \mathcal{F}^T\to \mathcal{F}^T$,
\begin{equation*}
(R^Tf)(t):=u^f(0,t)=\langle {\rm see\,}
(\ref{u_2_Duhamel})\rangle=\int_0^t r(t-s)f(s)\,ds, \quad
0\leqslant t\leqslant T,
\end{equation*}
where a \emph{response function} $r(t):=r(0,t)$ is smooth on
$[0,T]$ and satisfies $r(0)<0$ by (\ref{jump of r}). The operator
$R^T$ is compact.

An {\bf extended response operator} $\widetilde R^{2T}:
\mathcal{F}^{2T}\to \mathcal{F}^{2T}$ is introduced via the
problem $\widetilde 1^{2T}_{\rm Neum}$ by the rule
\begin{equation}
\label{R_2T_Neum} (\widetilde R^{2T}f)(t):=u^f(0,t)=\int_0^t
r(t-s)f(s)\,ds, \qquad 0\leqslant t\leqslant 2T;
\end{equation}
it is also a compact operator.

Comparing the response operators of the systems $\alpha^T_{\rm
Dir}$ and $\alpha^T_{\rm Neum}$, it is easy to see that $ R_{\rm
Neum}=(R_{\rm Dir})^{-1}$.
\smallskip

A {\bf connecting operator} $C^T:\mathcal{F}^T\to \mathcal{F}^T$
is $ C^T:=(W^T)^*W^T$. Along with $W^T,$ it is a compact operator.
\begin{theorem}
\label{C_T_theor_NM} The relation
\begin{equation*}
C^T=-\frac{1}{2}\,(S^T)^*\widetilde
R^{2T}J^{2T}S^T
\end{equation*}
and representation
\begin{equation}
\label{C_T_NM_repr1}
(C^Tf)(t)=-\int_0^T\left[\frac{1}{2}\int_{|t-s|}^{2T-t-s}r(\eta)\,d\eta\right]f(s)\,ds,\qquad
0\leqslant t\leqslant T
\end{equation}
are valid.
\end{theorem}
We omit the proof, which is quite analogous to the proof of
Theorem \ref{C_T_theor}.

In accordance with Lemma (\ref{Lemma_NM_contr}), the operator
\begin{equation*}
\left(W^T\frac{d}{dt}\right)^*\left(W^T\frac{d}{dt}\right)=\left(\frac{d}{dt}\right)^*(W^T)^*W^T
\frac{d}{dt}=-\frac{d}{dt}\,C^T \frac{d}{dt}
\end{equation*}
is a positive definite isomorphism in $\mathcal{F}^T$ \footnote{In
these equalities, $\left(\frac{d}{dt}\right)^*=-\frac{d}{dt}$ is
understood as the conjugate to $\frac{d}{dt}: \mathcal{F}^T\to
\mathcal{F}^T$, $\operatorname{Dom}\frac{d}{dt}=\left\{f\in
H^1[0,T]\,|\, f(0)=0\right\}$}. Simple calculations with regard to
(\ref{C_T_NM_repr1}) provide the representation
\begin{align} \label{C_T_dt_NM}
&\left(\left(-\frac{d}{dt}\,C^T\frac{d}{dt}\right)f\right)(t)=-r(0)f(t)\,+\notag \\
&\frac{1}{2}\int_0^T\left[\,r^\prime(2T-t-s)+r^\prime(|t-s|)\,\right]f(s)\,ds,
\quad 0\leqslant t\leqslant T.
\end{align}

\subsection{BCP for $\alpha^T_{\rm Neum}$.}
In both of the systems $\alpha^T_{\rm Dir}$ and $\alpha^T_{\rm
Neum}$, a control $f$ determines not only the wave
$u^f(\,\cdot\,,t)$ but its velocity $u^f_t(\,\cdot\,,t)$, i.e., a
\emph{complete dynamical state} that is a pair
$\{u^f(\,\cdot\,,t),u^f_t(\,\cdot\,,t)\}$. By injectivity of the
control operators $W^T_{\rm Dir}$ and $W^T_{\rm Neum}$, each
component of the final state
$\{u^f(\,\cdot\,,T),u^f_t(\,\cdot\,,T)\}$ determines the control
uniquely. However, taking into account the concrete properties of
these operators, to have a well-posed BCP it is natural to control
$u^f(\,\cdot\,,T)=W^Tf$ in $\alpha^T_{\rm Dir}$ but
$u^f_t(\,\cdot\,,T)=W^T\frac{d}{dt}f$ in $\alpha^T_{\rm Neum}$.

Namely, the following setup of the BCP is relevant for
$\alpha^T_{\rm Neum}$:

\noindent $\bullet$\,\,\,{\it given a target state $y\in
\mathcal{H}^T$, find a control $f\in \mathcal{F}^T$ such that
$u^f_t(\,\cdot\,,T)=y$}.

\noindent This problem is equivalent to the equation
$W^T\frac{d}{dt}f=y$ and, by Lemma \ref{Lemma_NM_contr}, is
uniquely solvable for all $y\in \mathcal{H}^T$.
\smallskip

Now, consider a \emph{special BCP} for $\alpha^T_{\rm Neum}$. Let
$y\in \mathcal{H}^T$ satisfy $-y^{\prime \prime}+qy=0$ as $x \in
(0,x(T))$.
\begin{theorem}
\label{Theor_NM_SP} The control
$f=\left(W^T\frac{d}{dt}\right)^{-1}y$ is a unique solution to the
equation
\begin{equation}
\label{C_T_eqn_NM}
-\frac{d}{dt}\,C^T\frac{d}{dt}=-y(0)1^T+y^\prime(0)(R^T)^*1^T
\end{equation}
in $\mathcal{F}^T$, where $1^T(t)=1.$
\end{theorem}
We omit the proof, which is quite analogous to the derivation of
(\ref{C_T_eqn}) but recommend the reader to recover the details.

With regard to (\ref{C_T_dt_NM}), we see that (\ref{C_T_eqn_NM})
is a second order Fredholm integral equation, which is uniquely
solvable because $-\frac{d}{dt}\,C^T\frac{d}{dt}$ is a (positive
definite) isomorphism in $\mathcal{F}^T$. Its solution $f$ is of
the class $C^\infty[0,T]$.

For the solution $f$ of the special BCP, the relation (\ref{repr
Wdf/dt}) implies
\begin{equation*}
y(x)=u^f_t(x,T)=r(x,\tau(x)+0)f(T-\tau(x))+\int_0^{T-\tau(x)}r_t(x,T-s)f(s)\,ds.
\end{equation*}
Putting here $x=x(T)$ and taking into account (\ref{jump of r}),
we get the equality
\begin{equation}
\label{f_NM_sol}
y(x(T))=-\left[\rho(0)\rho(x(T))\right]^{-\frac{1}{4}}f(0),
\end{equation}
which will be used in inverse problems.

The extended response operator $\widetilde R^{2T}$ of the system
$\alpha^T_{\rm Neum}$ determines the operator $\widetilde
R^{2\xi}$ and $R^{\xi}$ for all $\xi\leqslant T.$ Along with them,
$\widetilde R^{2T}$ determines the solutions $f=f^\xi(t)$ of the
family of the equations
\begin{equation}
\label{C_T_eqn_NM_xi}
-\frac{d}{dt}\,C^\xi\frac{d}{dt}=-y(0)1^\xi+y^\prime(0)(R^\xi)^*1^\xi
\end{equation}
in $\mathcal{F}^\xi$, the relation
\begin{equation}
\label{f_NM_sol_xi}
y(x(\xi))=-\left[\rho(0)\rho(x(\xi))\right]^{-\frac{1}{4}}f^\xi(0),\qquad
0\leqslant \xi\leqslant T
\end{equation}
being valid by (\ref{f_NM_sol}).

\section{Inverse problems}
\subsection{IP for $\alpha^T_{\rm Dir}$}
Consider the system $\alpha^T_{\rm Dir}$ with $\rho=1$ and a
variable $q$. The dynamical inverse problem is set up as follows:

\noindent $\bullet$\,\,\,{\it given the operator $\widetilde
R^{2T}$, recover the potential $q$ on the interval $[0,x(T)]$.}

\noindent Note that in this case we have $\tau(x)=x,$  $x(T)=T$,
whereas (\ref{R_2T_repr}) takes the form
\begin{equation}
\label{R_2T_repr_IP1} (\widetilde
R^{2T}f)(t)=-\frac{df}{dt}(t)+\int_0^t r(t-s)f(s)\,ds,\qquad
0\leqslant t\leqslant 2T.
\end{equation}
Hence, to give $\widetilde R^{2T}$ is to know the function $r(t)$,
$\,\,0\leqslant t\leqslant 2T$.
\smallskip

\noindent {\bf Family of equations}\,\,\, Let $y$ satisfies
\begin{align}
-&y''+qy=0, \qquad 0<x<x(\xi)\label{y_eq}\\
& y(0)=0,\,\, y'(0)=1. \label{y_cnd}
\end{align}
The equations (\ref{C_T_del_eqn}) take the form $C^\xi
f=\varkappa^\xi$ or
\begin{equation}
\label{C_T_del_eqn_GL}
f(t)+\int_0^\xi\left[p(2\xi-t-s)-p(|t-s|)\right]f(s)\,ds=\xi-t,\qquad
0\leqslant t\leqslant \xi,
\end{equation}
where $p(t):=\frac{1}{2}\int_0^t r(\eta)\,d\eta$ is known for
$t\in [0,2T].$ By (\ref{Jump_cnd_xi}), we have
\begin{equation}
\label{Jump_cnd_GL} y(\xi)=f^\xi(0),\qquad 0\leqslant \xi\leqslant
x(T)=T.
\end{equation}
\smallskip

\noindent {\bf Solving IP}\,\, Given the operator $\widetilde
R^{2T}$ (or, equivalently, the function $r|_{[0,2T]}$), one can
compose the equations (\ref{C_T_del_eqn_GL}), find their solutions
$f=f^\xi(t)$, and determine $y|_{[0,x(T)]}$ by
(\ref{Jump_cnd_GL}). Then one can recover the potential from
(\ref{y_eq}):
\begin{equation*}
q(\xi)=\frac{y''(\xi)}{y(\xi)},\qquad 0\leqslant \xi\leqslant
x(T)=T.
\end{equation*}
As is easy to see, possible presence of a finite number of zeros
of $y(\xi)$ does not prevent from recovering the smooth
$q|_{[0,T]}$. The IP is solved.

The conditions, which provide a solvability of this IP, are also
well-known: see, e.g., \cite{Bl71, Blag2, R}.
\begin{theorem}
\label{theor_BC_EQN} The operator defined by the r.h.s. of
$(\ref{R_2T_repr_IP1})$ is the extended response operator of a
system $\alpha^T_{\rm Dir}$ (with $\rho=1$) if and only if the
operator defined by the r.h.s. of $(\ref{C_T_repr1})$ (with
$\rho(0)=1$) is a positive definite isomorphism in
$\mathcal{F}^T$.
\end{theorem}
\smallskip

\noindent {\bf Gelfand-Levitan equation}\,\,\, Let $f=f^\xi(t)$
($0\leqslant t\leqslant \xi$) be the solutions to
(\ref{C_T_del_eqn_GL}). Denoting $f(\xi,t):=f^\xi(\xi-t)$ and
changing variables in (\ref{C_T_del_eqn_GL}), we get
\begin{equation*}
f(\xi,t)+\int_0^\xi\left[p(t+s)-p(|t-s|)\right]
f(\xi,s)\,ds=t,\quad 0\leqslant t\leqslant \xi.
\end{equation*}
Differentiating w.r.t. $\xi$, we have
\begin{equation*}
f_\xi(\xi,t)+\left[p(t+\xi)-p(|t-\xi|)\right]f(\xi,\xi)+\int_0^\xi\left[p(t+s)-p(|t-s|)\right]
f_\xi(\xi,s)\,ds=0.
\end{equation*}
Setting $L(\xi,t):=\left[f(\xi,\xi)\right]^{-1}f_\xi(\xi,t)$ and
$F(s,t):=p(t+s)-p(|t-s|)$, we arrive at the classical {\it
Gelfand-Levitan equation}
\begin{equation*}
F(\xi,t)+L(\xi,t)+\int_0^\xi
F(s,t)L(\xi,s)\,ds=0,\qquad 0\leqslant t\leqslant\xi
\end{equation*}
(see \cite{AM,GL}).
\smallskip

\noindent{\bf Transformation operator}\,\,\, In the original
Gelfand-Levitan approach, the function $L$ plays the role of the
kernel of a \emph{transformation operator} $V$ of the form
\begin{equation*}
(Vg)(x)=g(x)+\int_0^x L(x,t)g(t)\,dt,\qquad x\geqslant0.
\end{equation*}
This operator links solutions of the perturbed and unperturbed
spectral Sturm-Liouville problems
\begin{equation*}
\begin{cases}
&-\psi^{\prime \prime}+q\psi=\lambda\psi,\,\,\,x>0\\
&\psi(0)=0,\,\,\,\psi^\prime(0)=1
\end{cases}\quad {\rm and}\quad
\begin{cases}
&-\psi^{\prime \prime}_0=\lambda\psi_0,\quad x>0\\
&\psi_0(0)=0,\,\,\,\psi^\prime_0(0)=1
\end{cases}
\end{equation*}
by $\psi=V\psi_0$, i.e.
\begin{equation*}
\psi(x,\lambda)=\frac{\sin{\sqrt{\lambda}}x}{\sqrt{\lambda}}+\int_0^x
L(x,t)\,\frac{\sin{\sqrt{\lambda}}t}{\sqrt{\lambda}}\,dt, \quad
x\geqslant 0,
\end{equation*}
(see \cite{N}, Chapter VIII, sec. 26.5). Let us clarify the place
of the operator $V$ in the BC-method.

For the system $\alpha^T_{\rm Dir}$, consider the BCP of the form
$u^f(\,\cdot\,,T)=\psi(\,\cdot\,,\lambda),$ i.e.,
\begin{equation}
\label{BCP_trans} (W^Tf)(x)=\psi(x,\lambda),\qquad 0\leqslant
x\leqslant x(T)=T.
\end{equation}
By analogy with the proof of Theorem \ref{Theor_DIR_BCP}, one can
easily derive that the solution to (\ref{BCP_trans}) satisfies
\begin{equation*}
(C^Tf)(t)=\frac{\sin{\sqrt{\lambda}(T-t)}}{\sqrt{\lambda}},\qquad
0\leqslant t\leqslant T.
\end{equation*}
For $\lambda=0$, the latter equation coincides with
(\ref{C_T_eqn}). In the mean time, we have
\begin{equation*}
C^Tf=(W^T)^*W^Tf=\langle \,{\rm see\,}
(\ref{BCP_trans})\rangle=(W^T)^*\psi(\,\cdot\,,\lambda)=Y^T\psi_0(\,\cdot\,,\lambda),
\end{equation*}
where $Y^T: {\mathcal F}^T\to {\mathcal F}^T$,
$(Y^Tg)(t):=g(T-t),\,\,\, 0\leqslant t\leqslant T$. Hence,
$\psi(\,\cdot\,,\lambda)=\left[(W^T)^*\right]^{-1}Y^T\psi_0(\,\cdot\,,\lambda)$
holds in $\mathcal{H}$ and implies
\begin{equation*}
V=\left[(W^T)^*\right]^{-1}Y^T.
\end{equation*}
Thus, solving the BCP-equations (\ref{C_T_del_eqn_GL}), one can
find the function $L$ and determine the transformation operator
$V$.
\smallskip

\noindent{\bf Constructing $(W^T)^{-1}$}\,\, Assume that equations
(\ref{C_T_del_eqn_GL}) are solved. Their solutions $f^\xi$
determine the controls $f^{T,\,\xi}\in \mathcal{F}^T$, which solve
the BCPs (\ref{BCP_del}). Choosing a small enough $\Delta\xi>0$
and subtracting the controls, we obtain
\begin{equation*}
\left(W^T\left[\frac{f^{T,\,\xi}-f^{T,\,\xi-\Delta\xi}}{\Delta\xi}\right]\right)(x)=
\left\{\begin{array}l \frac{y(x)}{\Delta\xi}\,\,\,{\rm as}\,\,\,
x\in
[x(\xi-\Delta\xi),x(\xi)]\\
0\qquad\text{for other $x\in[0,x(T)]$}
\end{array}
\right. .
\end{equation*}
Since $x(\xi)=\xi$, the limit passage $\Delta \xi \to 0$ obviously
yields
\begin{equation}\label{**}
(W^T f_\xi^{T,\,\xi})(x)=y(\xi)\delta(x-\xi),\qquad x\in [0,x(T)],
\end{equation}
where the derivative $f_\xi^{T,\,\xi}$ has to be understood in the
sense of distributions. Taking into account a jump of
$f^{T,\,\xi}$ at $t=T-\xi$, we have
$$f_\xi^{T,\,\xi}(t)=f^{T,\,\xi}(T-\xi+0)\, \delta (t-(T-\xi))+
\widetilde f^{T,\,\xi}_\xi(t), \qquad 0\leqslant t \leqslant T,$$
where $\widetilde f^{T,\,\xi}_\xi$ is a classical derivative of
$f^{T,\,\xi}\in \mathcal{F}^{T,\,\xi}$ on $[0,T]$ for
$t\not=T-\xi$.

By (\ref{Jump_cnd_GL}), we have
$y(\xi)=f^\xi(0)=f^{T,\,\xi}(T-\xi+0)$ and can write (\ref{**}) in
the form
\begin{equation}
\label{BCP_delta} W^Tg^{T,\,\xi}=\delta^{x(\xi)} \qquad {\rm
on}\quad 0\leqslant x \leqslant x(T),
\end{equation}
where $\delta^{x(\xi)}$ is a Dirac measure supported at
$x=x(\xi)=\xi$, and
\begin{equation*}
g^{T,\,\xi}(t):=[f^{T,\,\xi}(T-\xi+0)]^{-1}f^{T,\,\xi}_\xi(t)=\delta(t-(T-\xi))+\widetilde
g^{T,\,\xi}(t), \quad 0\leqslant t \leqslant T
\end{equation*}
is a generalized (singular) control, $\widetilde g^{T,\,\xi}\in
\mathcal{F}^{T,\,\xi}$ being its regular part. Thus, $g^{T,\,\xi}$
solves a special BCP of the form (\ref{BCP_delta})\footnote{A BCP
with the Dirac measures as target functions is a cornerstone of
the BC-method: see \cite{BDAN87, B08, BV1}.}.

For a smooth $a\in \mathcal{H}^T$, one has
\begin{align*}
&\left(W^T\left[\int_0^Tg^{T,\,\xi}(\cdot)a(\xi)\,d\xi\right]\right)(x)=
\int_0^T(W^Tg^{T,\,\xi})(x)a(\xi)\,d\xi=\langle
{\rm see\,} (\ref{BCP_delta})\rangle=\\
&\int_0^{x(T)}\delta(x-\xi)a(\xi)\,d\xi=a(x),\quad 0\leqslant x
\leqslant x(T)=T.
\end{align*}
Hence,
\begin{equation*}
\int_0^T
g^{T,\,\xi}(t)a(\xi)\,d\xi=\left[(W^T)^{-1}a\right](t),\qquad
0\leqslant t\leqslant T.
\end{equation*}
Thus, constructing the controls $g^{T,\,\xi}(t)$ via the solutions
of (\ref{C_T_del_eqn_GL}), we determine the (generalized) kernel
of the operator $(W^T)^{-1}.$
\smallskip

\noindent{\bf Visualization of waves}\,\,The BC-method provides
one more option that we call a \emph{visualization}.

On smooth controls, define a \emph{visualizing functional}
$\delta^{T,\,\xi}$ by the (properly understood) relations
\begin{align*}
&\langle \delta^{T,\,\xi},f\rangle:=\left(C^T
g^{T,\,\xi},f\right)_{\mathcal{F}^T}=\left(W^Tg^{T,\,\xi},W^Tf\right)_{\mathcal{H}^T}=\\
&\langle {\rm see\,} (\ref{BCP_delta})\rangle=
\int_0^{x(T)}\delta(x-\xi)u^f(x,T)\,dx=u^f(\xi,T).
\end{align*}
Thus, constructing such functionals, we can recover (visualize)
the waves.

One more interpretation is the following. Expanding a control $f$
over the 'standard' Dirac measures
$\{\delta^\xi\}_{0\leqslant\xi\leqslant T}$, we get
$\int_0^T\delta(t-\xi)f(t)\,dt=f(\xi)$, i.e., nothing new. In the
mean time, expanding over
$\{\delta^{T,\,\xi}\}_{0\leqslant\xi\leqslant
T},\,\,\delta^{T,\,\xi}:=C^T g^{T,\,\xi}$, we get
$\int_0^T\delta^{T,\,\xi}(t)f(t)\,dt=u^f(\xi,T)$, i.e., {\it see
the wave} $u^f$. So, to visualize waves is to pass from
$\{\delta^\xi\}$ to $\{\delta^{T,\,\xi}\}$, i.e., choose a proper
continual basis in the outer space $\mathcal{F}^T$.

\subsection{IP for $\alpha^T_{\rm Neum}$}
Consider the system $\alpha^T_{\rm Neum}$ with a variable $\rho$
and $q=0$. Recall the assumption
$\int_0^\infty\rho^{\frac{1}{2}}(x)\,dx=\infty$, which provides
the well-posedness of the problem $\widetilde 1^{2T}_{\rm Neum}$
for all $T>0$.

The inverse problem is set up as follows:

\noindent $\bullet$\,\,\,{\it given the operator $\widetilde
R^{2T}$,  recover the density $\rho$ on the interval $[0,x(T)]$.}

\noindent Note that to give $\widetilde R^{2T}$ is to know the
function $r(t)$, $0\leqslant t\leqslant 2T$ (see
(\ref{R_2T_Neum})).
\smallskip

\noindent {\bf Krein equation}\,\,\, Any linear function $y=ax+b$
satisfies $-y^{\prime \prime}=0$ and, hence, can be used as a
target function in the special BCP for $\alpha^T_{Neum}$ (see sec.
3.4). Put $y=-1.$ In this case, with regard to (\ref{C_T_dt_NM})
the family (\ref{C_T_eqn_NM_xi}) takes the form
\begin{equation}
\label{C_T_Krein_NM}
-r(0)f(t)+\frac{1}{2}\int_0^\xi\left[r^\prime(2\xi-t-s)+r^\prime(|t-s|)\right]f(s)\,ds=1,\,\,0\leqslant
t\leqslant \xi.
\end{equation}
Extend the solution and the r.h.s. of (\ref{C_T_Krein_NM}) from
$0\leqslant t\leqslant \xi$ to $\xi\leqslant t\leqslant 2\xi$ by
evenness with respect to $t=\xi$ and then substitute $t\mapsto
t-\xi$. One can easily check that the function
\begin{equation}
\label{g_def} g=g(\xi,t)\,:=\,\left\{\begin{array}l f(\xi+t)\quad
-\xi\leqslant t\leqslant
0\\
f(\xi-t)\quad 0\leqslant t\leqslant \xi
\end{array}
\right.
\end{equation}
satisfies
\begin{equation*}
-r(0)g(t)+\frac{1}{2}\int_{-\xi}^\xi
r'(|t-s|)g(s)\,ds=1,\qquad -\xi\leqslant t\leqslant \xi
\end{equation*}
that is the classical {\it M.Krein equation}. It was derived for
solving the inverse spectral problem for inhomogeneous string in
\cite{Kr1,Kr2}, see also \cite{Bl71} and \cite{GoSh1,GoSh2}.
\smallskip

\noindent {\bf Solving IP}\,\,\, Return to equations
(\ref{C_T_Krein_NM}) and find their solutions $f=f^\xi(t)$. By
(\ref{f_NM_sol_xi}), we have
$f^\xi(0)=-\left[\rho(0)\rho(x(\xi))\right]^{\frac{1}{4}}$. Since
$\rho(0)=r^2(0)$ is known, we determine
\begin{equation*}
\rho(x(\xi))=\left[r(0)\right]^{-1}\left[f^\xi(0)\right]^4,\qquad
0\leqslant \xi\leqslant T.
\end{equation*}
Then, we can find the function
\begin{equation*}
x(\xi)=\int_0^\xi\rho^{-\frac{1}{2}}(x(s))\,ds,\qquad
0\leqslant\xi\leqslant T
\end{equation*}
and its inverse function (eikonal) $\tau(x)$ on $0\leqslant
x\leqslant x(T)$. At last, we recover
\begin{equation*}
\rho(x)=\tau^\prime(x),\qquad 0\leqslant x\leqslant x(T).
\end{equation*}
The IP is solved.

The solvability conditions for this IP are also well known: see
\cite{Bl71,R}.
\begin{theorem}
\label{Theor_IP_Neum} The operator defined by the r.h.s. of
$(\ref{R_2T_Neum})$ is the extended response operator of a system
$\alpha^T_{\rm Neum}$ (with $q=0$) if and only if the operator
defined by the r.h.s. of $(\ref{C_T_dt_NM})$ is a positive
definite isomorphism in $\mathcal{F}^T$.
\end{theorem}
By analogy with $\alpha^T_{\rm Dir}$, the solutions of
(\ref{C_T_Krein_NM}) can be used for constructing $(W^T)^{-1}$ and
visualization of waves in $\alpha^T_{\rm Neum}$.
\smallskip

\noindent{\bf Pariiskii equation}\,\,\,In the special BCP for
$\alpha^T_{\rm Neum}$, put $y=-x$ as a target function. The
equations (\ref{C_T_eqn_NM_xi}) take the form
\begin{align}
& -r(0)f(t)+\frac{1}{2}\int_0^\xi\left[r^\prime(2\xi-t-s)+r^\prime(|t-s|)\right]f(s)\,ds\notag\\
& =-\int_t^\xi r(s-t)\,ds,\qquad 0\leqslant t\leqslant
\xi.\label{Parii_eq}
\end{align}
Arguing as in the derivation of the Krein equation, we pass from
(\ref{Parii_eq}) to
\begin{equation*}
-r(0)g(t)+\frac{1}{2}\int_{-\xi}^\xi
r^\prime(|t-s|)g(s)\,ds=-\int_0^{|t|} r(s)\,ds,\qquad
-\xi\leqslant t\leqslant \xi,
\end{equation*}
where $g=g(\xi,t)$ is defined as in (\ref{g_def}). Differentiating
w.r.t. $\xi$, with regard to the evenness of $g$, we get
\begin{equation*}
-r(0)g_\xi(\xi,t)+\frac{1}{2}\int_{-\xi}^\xi r'(|t-s|)
g_\xi(\xi,s)\,ds+\frac{g(\xi,\xi)}{2}\left[r^\prime(\xi+t)+r^\prime(\xi-t)\right]=0.
\end{equation*}
Then, denoting $G(\xi,t):=[g(\xi,\xi)]^{-1} g_\xi(\xi,t)$, we
arrive at the Pariiskii equation (see \cite{Bl71}):
\begin{equation*}
\label{Par_eqn} -r(0)G(t)+\frac{1}{2}\int_{-\xi}^\xi
r^\prime(|t-s|)
G(s)\,ds+\frac{1}{2}\left[r^\prime(\xi+t)+r^\prime(\xi-t)\right]=0,\qquad
-\xi\leqslant t\leqslant \xi.
\end{equation*}

\section{Inverse scattering problem}
\subsection{System $\alpha_{\rm scatt}$} In the rest of the paper,
Convention 1 is cancelled.

A dynamical system $\alpha_{\rm scatt}$ is associated with the
(forward) problem
\begin{align}
& u_{tt}-u_{xx}+qu=0, && x>0,\,\,\, -\infty<t<x\label{sc_eq}\\
& u|_{t<-x}=0\label{sc_c1}\\
& \lim_{s\to\infty}u(s,\tau-s)=f(\tau), && \tau\geqslant
0,\label{sc_c2}
\end{align}
where $q\in C^\infty[0,\infty)$ is a \emph{potential}, $f$ is a
\emph{control}, $u=u^f(x,t)$ is a solution (\emph{wave}). For the
sake of simplicity, we assume that the potential is compactly
supported: $\operatorname{supp}q\subset [0,a]$, $a<\infty$. Under
this assumption, (\ref{sc_c2}) implies $u^f|_{t<-a}=f(x+t)$. By
this, a control can be interpreted as an incident wave incoming
from $x=\infty.$

The scattering problem (\ref{sc_eq})--(\ref{sc_c2}) is well-posed:
this is established by reducing to relevant Volterra-type integral
equations (see, e.g., \cite{Blag2}). Note that, owing to
hyperbolicity, this setup does not need to impose a boundary
condition at $x=0$. In a sense, it is a natural analog of the
extended problems $\widetilde 1^{2T}_{\rm Dir}$ and $\widetilde
1^{2T}_{\rm Neum}$, which need no condition at $x=x(T)$. If $f\in
C^\infty[0,\infty)$ vanishes near $\tau=0$, the problem has a
unique classical smooth solution $u^f$.

Since $q|_{x>a}=0$, for large $x$'s the solution satisfies
$u_{tt}-u_{xx}=0$ and, hence, is a sum of two D'Alembert waves:
\begin{equation}
\label{Sc_sol_repr} u^f(x,t)|_{x>a}=f(x+t)+f^*(x-t),
\end{equation}
where the second summand describes the wave reflected by the
potential and outgoing to $x=\infty.$

Taking $f=\delta(t)$, one can introduce a fundamental solution of
the form $u^\delta(x,t)=\delta(t+x)+w(x,t)$, which satisfies
\begin{equation*}
u^\delta(x,t)|_{x>a}=\delta(x+t)+r(x-t)
\end{equation*}
with $r\in C^\infty[0,\infty),$
$\operatorname{supp\,}r\subset[0,2a]$. The Duhamel representation
$u^f=u^\delta*f$ holds for the classical solutions. Also, it
provides the relevant definition of a generalized solution for $f
\in L_2(0,\infty)$ and leads to
\begin{equation}
\label{Sc_Duhamel} f^*(\tau)=\int_0^\infty
r(\tau+s)f(s)\,ds,\qquad \tau\geqslant 0.
\end{equation}
Note that $\operatorname{supp}f^*\subset [0,2a],$ so that the
reflected wave $f^*(x-t)$ in (\ref{Sc_sol_repr}) is compactly
supported on $t\leqslant x \leqslant \infty$ for any $t$. Also, we
have
\begin{equation}
\label{Sc_sol_repr2} u^f(x,t)|_{t<-a}=f(x+t),
\end{equation}
which means that for these times the incident wave does not
interact with the potential and no reflected waves are produced.
\smallskip

\noindent An {\bf outer space} of the system $\alpha_{\rm scatt}$
is the space of controls $\mathcal{F}:=L_2(0,\infty)$. It contains
the subspaces $\mathcal{F}^\xi=\left\{f\in
\mathcal{F}\,|\,\operatorname{supp}f\subset [\xi,\infty)\right\}$
($\xi\geqslant 0$) formed by the delayed controls.
\smallskip

\noindent An {\bf inner space} is $\mathcal{H}:=L_2(0,\infty)$ (of
functions of $x$). It contains the subspaces
$\mathcal{H}^\xi:=\left\{y\in
\mathcal{H}\,|\,\operatorname{supp}y\subset
[\xi,\infty)\right\},\,\,\,\,\xi\geqslant 0$.
\smallskip

\noindent A {\bf control operator} $W:\mathcal{F}\to \mathcal{H}$
acts by the rule
\begin{equation*}
(Wf)(x):=u^f(x,0),\qquad x\geqslant 0.
\end{equation*}
It maps $\mathcal{F}$ onto $\mathcal{H}$ isomorphically, the
relation $W \mathcal{F}^\xi=\mathcal{H}^\xi$ being valid for all
$\xi\geqslant 0.$ These facts are derived from the representation
\begin{equation}
\label{Sc_W_repr} (Wf)(x)=f(x)+\int_0^x w(x,-s)f(s)\,ds,\qquad
x\geqslant 0,
\end{equation}
which, in turn, follows from the Duhamel representation.
\smallskip

\noindent A {\bf connecting operator}
$C:\mathcal{F}\to\mathcal{F},$
\begin{equation*}
C:=W^*W
\end{equation*}
is a positive definite isomorphism in $\mathcal{F}$. It connects
the metrics of the outer and inner spaces:
\begin{equation}
\label{Sc_C}
(Cf,g)_{\mathcal{F}}=(Wf,Wg)_\mathcal{H}=\left(u^f(\,\cdot\,,0),u^g(\,\cdot\,,0)\right)_\mathcal{H}.
\end{equation}
\smallskip

\noindent A {\bf response operator} of the system $\alpha_{\rm
scatt}$ is $R:\mathcal{F}\to\mathcal{F},$
\begin{equation*}
(Rf)(\tau):=\lim_{s\to
+\infty}u^f(s,s-\tau),\qquad \tau\geqslant 0.
\end{equation*}
For $f\in \mathcal{F}$ vanishing at $\infty$, by
(\ref{Sc_sol_repr}), this limit is $f^*(\tau).$ Hence, by
(\ref{Sc_Duhamel}) we get
\begin{equation}
\label{Sc_R_repr} (Rf)(\tau)=\int_0^\infty
r(\tau+s)f(s)\,ds,\qquad \tau\geqslant 0.
\end{equation}
Thus, $R$ is a compact self-adjoint operator.
\smallskip

\noindent {\bf Basic relation}\,\, As ever in the BC-method, the
key fact is an explicit formula, which relates the connecting and
response operators.
\begin{theorem}
\label{Sc_teor_C_repr} The equality
\begin{equation}
\label{Sc_C_repr} C=\mathbb{I}+R
\end{equation}
holds, where $\mathbb{I}$ is the identity operator in
$\mathcal{F}$.
\end{theorem}
\begin{proof}\,\, Fix two controls $f,g\in
C^\infty_0(0,\infty)$. Analyzing the positions of
$\operatorname{supp}u^f(\,\cdot\,,s)$ and
$\operatorname{supp}u^g(\,\cdot\,,t)$ at the semi-axis $x
\geqslant 0$, one can easily conclude that the Blagoveschenskii
function
\begin{equation*}
b(s,t):=\int_0^\infty u^f(x,s)u^g(x,t)\,dx
\end{equation*}
is well defined if (and only if) $s+t\leqslant 0$ holds. For the
admissible $s,t$ it satisfies
\begin{align}
& b_{tt}(s,t)-b_{ss}(s,t)=\int_0^\infty\left[u^f(x,s)u^g_{tt}(x,t)-u^f_{ss}(x,s)u^g(x,t)\right]\,dx=\notag\\
&\int_0^\infty\left[u^f(x,s)\left(u^g_{xx}(x,t)-q(x)u^g(x,t)\right)-\left(u^f_{ss}(x,s)-
q(x)u^f(x,s)\right)u^g(x,t)\right]\,dx=\notag\\
&\left[u^f(x,s)u_x^g(x,t)-u^f_x(x,s)u^g(x,t)\right]\bigr|_{x=0}^{x=\infty}=0,\qquad
s+t<0\label{Sc_b_cnd0}
\end{align}
by compactness of the supports of the waves.

Taking into account (\ref{Sc_sol_repr2}), we have
\begin{align*}
& b(t,t)|_{t<-a}=\int_0^\infty u^f(x,t)u^g(x,t)\,dx=\int_{-t}^\infty f(x+t)g(x+t)\,dx=\\
& \int_0^\infty f(\tau)g(\tau)\,d\tau=(f,g)_\mathcal{F}.
\end{align*}
Therefore,
\begin{equation}
\label{Sc_b_cnd} \lim_{t\to-\infty}b(t,t)=(f,g)_\mathcal{F}.
\end{equation}
In accordance with (\ref{Sc_sol_repr}) and (\ref{Sc_sol_repr2}),
we have
\begin{align*}
& \lim_{t\to-\infty}b(-t,t)=\lim_{t\to-\infty}\int_0^\infty
u^f(x,-t)u^g(x,t)\,dx=\notag\\
&\lim_{t\to-\infty}\int_{-t}^\infty f^*(x+t)g(x+t)\,dx=
\int_0^\infty f^*(\tau)g(\tau)\,d\tau=(Rf,g)_\mathcal{F}
\label{Sc_b_cnd1}
\end{align*}
and, quite analogously,
\begin{equation}\label{888}
\lim_{t\to-\infty}b(t,-t)=(f,Rg)_\mathcal{F}=(Rf,g)_\mathcal{F}
\end{equation}
in view of $R^*=R.$ Also, the positions of
$\operatorname{supp}u^f(\,\cdot\,,s)$ and
$\operatorname{supp}u^g(\,\cdot\,,t)$ in $x\geqslant 0$ easily
imply
\begin{equation}
\label{Sc_b_cnd2}
\lim_{s\to-\infty}b(s,0)=\lim_{s\to-\infty}\int_0^\infty
u^f(x,s)u^g(x,0)\,dx=0.
\end{equation}
Since the function $b$ satisfies $b_{tt}-b_{ss}=0$ (see
(\ref{Sc_b_cnd0})), it is a sum of the D'Alembert solutions:
\begin{equation*}
b(s,t)=\Phi(t+s)+\Psi(t-s),\quad t+s\leqslant 0
\end{equation*}
with $\Phi(\tau)$ and $\Psi(\tau)$ defined on $\tau\leqslant 0$
and $-\infty<\tau<\infty$ respectively.

The relations (\ref{Sc_b_cnd}), (\ref{888}), (\ref{Sc_b_cnd2})
imply
\begin{equation*}
\Phi(-\infty)+\Psi(0)=(f,g)_\mathcal{F},\quad
\Phi(0)+\Psi(\infty)=(Rf,g)_\mathcal{F},\quad
\Phi(-\infty)+\Psi(\infty)=0
\end{equation*}
respectively. As a result, we have
\begin{eqnarray*}
\left((\mathbb{I}+R)f,g\right)_\mathcal{F}=(f,g)_\mathcal{F}+(Rf,g)_\mathcal{F}=
\Phi(-\infty)+\Psi(0)+\Phi(0)+\Psi(\infty)\\
=\Psi(0)+\Phi(0)=b(0,0)=\left(u^f(\,\cdot\,,0),u^g(\,\cdot\,,0)\right)_\mathcal{F}=\langle
{\rm see\,} (\ref{Sc_C}) \rangle=(C f,g)_\mathcal{F}.
\end{eqnarray*}
By arbitrariness of the choice of $f$ and $g$, we arrive at
$\mathbb{I}+R=C.$
\end{proof}

\subsection{Control problem} A natural setup of a control problem (CP) for the
system $\alpha_{\rm scatt}$ is the following:

\noindent $\bullet$\,\,{\it given $y\in \mathcal{H}$, find $f\in
\mathcal{F}$ such that $u^f(\,\cdot\,,0)=y$.}

\noindent This problem is equivalent the equation $Wf=y$ and is
well posed by the properties of the control operator $W$ (see
(\ref{Sc_W_repr})): for any $y\in \mathcal{H}$, it has a unique
solution $f=W^{-1}y\in \mathcal{F}$.

Now, let us consider a \emph{special CP}: take $y$, which
satisfies
\begin{align}
 -&y''+qy=-k^2y,\qquad x>0\quad (k>0)\label{Sct_y}\\
& y|_{x>a}=e^{-kx}. \label{Sct_y1}
\end{align}
\begin{theorem}
\label{Sc_control_teor} The control $f=W^{-1}y,$ which solves the
special CP, is the solution of the equation $Cf=e^{-k(\cdot)}$ in
$\mathcal{F}$, i.e., satisfies
\begin{equation}
\label{Sc_C_eq} f(\tau)+\int_0^\infty
r(\tau+s)f(s)\,ds=e^{-k\tau},\qquad \tau\geqslant 0.
\end{equation}
\end{theorem}
\begin{proof}
Choose a control $g\in C_0^\infty(0,\infty)$ (so that $g$ vanishes
near $t=0$) and denote $A:=\sup\operatorname{supp}g$,
$B:=\sup\operatorname{supp} u^g(\,\cdot\,,0)$. Note that for a
fixed $x>0$, the solution $u^g(x,t)$ vanishes as $t\leqslant-x$.
By the latter, we can represent
\begin{eqnarray}
u^g(x,0)=-\int_{-x}^0\frac{{\rm sh}\,
kt}{k}\left[u^g_{tt}(x,t)-k^2u^g(x,t)\right]\,dt=\notag\\
-\int_{-x}^0\frac{{\rm sh}\,
kt}{k}\left[u^g_{xx}(x,t)-(k^2+q(x))u^g(x,t)\right]\,dt.\label{Sc_u_o}
\end{eqnarray}
Fix an $N>\max{\{A,B,a\}}$. For the solution $f$ of the special
CP, one has
\begin{align*}
& (C f,g)_\mathcal{F}=\langle {\rm see\,}
(\ref{Sc_C})\rangle=\left(u^f(\cdot,0),u^g(\cdot,0)\right)_\mathcal{H}=
\left(y,u^g(\cdot,0)\right)_\mathcal{H}=\int_0^N y(x)u^g(x,0)\,dx=\\
&\langle {\rm see\,} (\ref{Sc_u_o})\rangle=-\int_0^N
dx\,y(x)\int_{-x}^0\frac{{\rm sh}\,
kt}{k}\left[u^g_{xx}(x,t)-(k^2+q(x))u^g(x,t)\right]\,dt=\\
&-\int_{-N}^0 dt\,\frac{{\rm sh}\, kt}{k}\int_{-t}^N
y(x)\left[u^g_{xx}(x,t)-(k^2+q(x))u^g(x,t)\right]\,dx=\\
&-\int_{-N}^0 dt\,\frac{{\rm sh}\,
kt}{k}\Bigl\{y(N)u_x^g(N,t)-y'(N)u^g(N,T)\Bigr.-\\
&\Bigl.-\int_{-x}^0\left[y''(x)-(k^2+q(x))y(x)\right]u^g(x,t\,dx)\Bigr\}=\\
&\langle {\rm see\,} (\ref{Sct_y}),
(\ref{Sct_y1})\rangle=-\int_0^N \frac{{\rm sh}\,
kt}{k}\left[e^{-kN}u^g_x(N,t)+ke^{-kN}u^g(N,t)\right]\,dt
\end{align*}
Increasing (if necessary) $N$, we can provide $u^g(N,t)=g(N+t)$
and $u^g_x(N,t)=g^\prime(N+t).$ For such $N$, we have
\begin{align*}
&(Cf,g)_\mathcal{F}=-\int_{-N}^0 \frac{{\rm sh}\,
kt}{k}\left[e^{-kN}g'(N+t)+ke^{-kN}g(N+t)\right]\,dt=\\
&-\int_{-N}^0 \left[({\rm ch}\,{kt})\,e^{-kN}-({\rm
sh}\,{kt})\,e^{-kN}\right]g(N+t)\,dt=\int_0^N
e^{-k\tau}g(\tau)\,d\tau=(e^{-k(\cdot)},g)_\mathcal{F}.
\end{align*}
By arbitrariness of $g$, we arrive at $Cf=e^{-k(\cdot)}$.
\end{proof}

{\bf Family of equations}\,\,\, Let $y$ be the solution of
(\ref{Sct_y}), (\ref{Sct_y1}). For a $\xi>0$, define a cut-off
function $y^\xi \in \mathcal{H}^\xi$ by
\begin{equation*}
y^\xi(x):=\left\{\begin{array}l
0,\qquad 0\leqslant x<\xi\\
y(x),\quad \xi\leqslant x<\infty
\end{array}.
\right.
\end{equation*}
With regard to the properties of the control operator $W$, the CP
$u^f(\,\cdot\,,0)=y^\xi$ is uniquely solvable, its solution
$f=f^\xi$ belonging to the subspace $\mathcal{F}^\xi$. By perfect
analogy with Theorem \ref{Sc_control_teor}, one can derive that it
satisfies the equation
\begin{equation}
\label{Sc_C_eq_del} f^\xi(\tau)+\int_\xi^\infty
r(\tau+s)f^\xi(s)\,ds=e^{-k\tau},\qquad \xi\leqslant\tau<\infty.
\end{equation}
On the other hand, representing by (\ref{Sc_W_repr})
\begin{equation*}
\left(Wf^\xi\right)(x)=f^\xi(x)+\int_\xi^x
w(x,-s)\,f^\xi(s)\,ds=y^\xi(x),\qquad x\geqslant 0
\end{equation*}
and putting $x=\xi$, we see that
\begin{equation}
\label{Sc_jump} f^\xi(\xi+0)=y(\xi),\qquad \xi \geqslant 0.
\end{equation}

\subsection{Inverse problem.} The dynamical inverse scattering
problem for the system $\alpha_{\rm scatt}$ is set up as follows:

\noindent$\bullet$\,\,{\it given the response operator $R$
(equivalently, the function $r|_{\tau\geqslant 0}$), recover the
potential $q|_{x\geqslant 0}$.}

To solve this IP, one can find the solutions $f^\xi$ of equations
(\ref{Sc_C_eq_del}) for all $\xi\geqslant 0,$ determine the
function $y$ by (\ref{Sc_jump}) and recover the potential by
\begin{equation*}
q(\xi)=\frac{y''(\xi)}{y(\xi)}-k^2,\qquad \xi \geqslant 0.
\end{equation*}
The solvability conditions for the IP are also well-known.
\begin{theorem}$(V.A.Marchenko)$\,
The operator defined by the r.h.s. of $(\ref{Sc_R_repr})$ is the
response operator of a system $\alpha_{\rm scatt}$ with a smooth
compactly supported potential $q$ if and only if the function $r$
is smooth and compactly supported on $\tau\geqslant 0$ and the
operator defined by the r.h.s. of $(\ref{Sc_C_repr})$ is a
positive definite isomorphism in $\mathcal{F}$.
\end{theorem}
\smallskip

\noindent {\bf Marchenko equation}\,\, Differentiating the
function $f(\xi,\tau):=f^\xi(\tau),\,\,0\leqslant t\leqslant \xi <
\infty$ w.r.t. $\xi$ in (\ref{Sc_C_eq_del}), we obtain
\begin{equation*}
f_\xi(\xi,\tau)-r(\tau+\xi)f^\xi(\xi)+\int_\xi^\infty
r(\tau+s)\,f_\xi(\xi, s)\,ds=0,\qquad \xi\leqslant\tau<\infty.
\end{equation*}
Denoting $g(\tau,\xi):=[f^\xi(\xi)]^{-1}f_\xi(\xi,\tau)$, we
arrive at
\begin{equation*}
\label{Sct_Mar} g(\xi,\tau)+\int_\xi^\infty r(\tau+s)
g(\xi,s)\,ds=r(\tau+\xi),\qquad \xi\leqslant\tau<\infty
\end{equation*}
that is the classical Marchenko equation: see, e.g., \cite{AgrM}.
\smallskip

\noindent{\bf Comments}

\noindent $\bullet$\,Analyzing the procedure, which recovers the
potential, one can see that to determine $q|_{x \geqslant\, \xi}$
it suffices to know $r|_{\tau \geqslant\, \xi}$ (so that to know
$r$ {\it for all} $\tau \geqslant 0$ is not necessary). Such a
{\it locality} is an inherent feature and advantage of the
BC-method. It is in consent with the fundamental physical property
of dynamical systems under consideration, which is a finiteness of
the wave propagation speed.

\noindent $\bullet$\, Writing (\ref{Sc_C_eq}) in the form
$W^*Wf=e^{-k(\cdot)},$ with regard to $Wf=y,$ we have
$y=\left(W^*\right)^{-1}e^{-k(\cdot)}$. Hence,
$\left(W^*\right)^{-1}$ is a {\it transformation operator}, which
maps the solution $e^{-k(\cdot)}$ of the unperturbed problem
(\ref{Sct_y}), (\ref{Sct_y1}) (with $q=0$) to the solution $y$ of
the perturbed problem.

\noindent $\bullet$\, Using the solutions $f^\xi$ of
(\ref{Sc_C_eq_del}) in the same way as at the end of sec. 4.1, one
can construct the singular controls $\delta^{0,\,\xi}$ (analogs of
$\delta^{T,\,\xi}$), which solve the CP of the form $W
\delta^{0,\,\xi}=\delta^\xi$. Then, the continual basis
$\{\delta^{0,\,\xi}\}_{0 \leqslant \xi <\infty}$ can be used for
visualization of waves.


\begin{thebibliography}{}

\bibitem{AgrM}
Z.S.Agranovich, V.A.Marchenko.
\newblock{The inverse problem of scattering theory.}
\newblock{\em New York, London: Gordon and Breach}, 1963.

\bibitem{AM}
S. A. Avdonin, V. S. Mikhaylov.
\newblock{The boundary control approach to inverse spectral theory.}
\newblock{\em Inverse Problems}, 26, no. 4, 045009, 19 pp. 2010.

\bibitem{BDAN87}
M.I.Belishev.
\newblock {On an approach to multidimensional inverse problems for
the wave equation.}
\newblock {\em Soviet Mathematics. Doklady}, 36, no 3, 481--484, 1988.

\bibitem{B87}
M.I.Belishev.
\newblock {The Gelfand--Levitan type equations in multidimensional
inverse problem for the wave equation.}
\newblock {\em Zapiski Nauch. Semin. LOMI}, 165, 15--20, 1987
(in Russian); English translation: \newblock{\em J. Sov. Math.,}
50, no 6, 1990.

\bibitem{B97}
M. I. Belishev.
\newblock{ Boundary control in reconstruction of manifolds and metrics
(the BC method)}
\newblock{\em  Inverse Problems}, 13, no 5, R1--R45, 1997.

\bibitem{B07}
M. I. Belishev.
\newblock{Recent progress in the boundary control
method.}
\newblock{\em Inverse Problems}, 23, no 5, R1--R67, 2007.

\bibitem{B08}
M. I. Belishev.
\newblock{Boundary control and inverse
problems: a one-dimensional version of the boundary control
method.}
\newblock{\em Zap. Nauchn. Sem. S.-Peterburg. Otdel. Mat.
Inst. Steklov. (POMI),} 354, 19--80, 2008 (in Russian); English
translation: \newblock{\em J. Math. Sci. (N. Y.),} 155, no 3,
343--378, 2008.

\bibitem{BV1}
M. I. Belishev, A. F. Vakulenko.
\newblock{Inverse problems on
graphs: recovering the tree of strings by the BC-method.}
\newblock{\em Journal of Inverse and Ill-Posed Problems,} 14, no
1, 29--46, 2006.

\bibitem{Bl71}
A. S. Blagoveschenskii.
\newblock{On a local approach to the
solution of the dynamical inverse problem for an inhomogeneous
string.}
\newblock{\em Trudy MIAN,} 115, 28--38, 1971 (in Russian).

\bibitem{Blag2}
A.S.Blagovestchenskii.
\newblock{Inverse Problems of Wave Processes.}
\newblock{\em VSP, Netherlands}, 2001.

\bibitem{GL}
I. M. Gel'fand, B. M. Levitan.
\newblock{On the determination
of a differential equation from its spectral function.}
\newblock{\em Izvestiya Akad. Nauk SSSR. Ser. Mat.,} 15, 309--360,
1951 (in Russian); English translation: \newblock{\em Amer. Math.
Soc. Transl. (2),} 1, 253--304, 1955.

\bibitem{GoSh1}
Gopinath B., Sondhi M. M.
\newblock{Determination of the shape of
the human vocal tract from acoustical measurements.}
\newblock{\em Bell Syst. Tech. J.,} 1195--1214, 1970.

\bibitem{GoSh2}
Gopinath B., Sondhi M. M.
\newblock{Inversion of the Telegraph
Equation and the Synthesis of Nonuniform Lines.}
\newblock{\em Proceedings of the IEEE,} 59, no 3, 383--392, 1971.

\bibitem{Kr1}
M. G. Krein.
\newblock{On the transfer function of a
one-dimensional boundary problem of the second order.}
\newblock{\em Dokl. Akad. Nauk. SSSR (N.S.),} 88, 405--408, 1953 (in
Russian).

\bibitem{Kr2}
M. G. Krein.
\newblock{On a method of efficient solution of
an inverse boundary problem.}
\newblock{\em Dokl. Akad. Nauk. SSSR
(N.S.),} 94, 987--990, 1954 (in Russian).

\bibitem{N}
M. A. Naimark.
\newblock{Linear Differentilal Operators.}
\newblock{\em Nauka, Moscow,} 1969 (in Russian).

\bibitem{R}
V. G. Romanov.
\newblock{Inverse problems of mathematical
physics.}
\newblock{\em Nauka, Moscow,} 1984. (in
Russian).
\end{thebibliography}
\end{document}